\documentclass{esannV2}
\usepackage[dvips]{graphicx}
\usepackage[latin1]{inputenc}
\usepackage{amssymb,amsmath,array}
\newtheorem{Lem}{\underline {Lemma}}
\newtheorem{Prop}{\underline {Proposition}}

\begin{document}
\title{Efficient Estimation of Multidimensional Regression Model with Multilayer Perceptron}

\author{Joseph Rynkiewicz $^1$
%
%
\vspace{.3cm}\\
%
Université Paris I - SAMOS/MATISSE\\
72 rue Regnault, Paris - France\\ 
}

\maketitle

\begin{abstract}
This work concerns estimation of multidimensional nonlinear regression
models using multilayer perceptron (MLP). The main problem with such model is that we have to know the covariance matrix of the noise to get optimal estimator. however we show that, if we choose as cost function  the logarithm of the determinant of the empirical error covariance matrix, we get an asymptotically optimal estimator. 
\end{abstract}

\section{Introduction}
Let us consider a sequence $\left(Y_{t},Z_{t}\right)_{t\in \mathbb{N}}$
of i.i.d.\footnote{It is not hard to extend all what we show in this paper for stationary mixing variables and so for time series} random vectors (i.e. identically distributed and independents). So, each couple $\left(Y_{t},Z_{t}\right)$ has the same law that a generic variable $(Y,Z)$. Moreover, we assume that the model can be written 
\[
Y_{t}=F_{W^0}(Z_{t})+\varepsilon _{t}
\]
where 

\begin{itemize}
\item $F_{W^0}$ is a function represented by a MLP with parameters or
weights $W^0$.
\item $(\varepsilon _{t})$ is an i.i.d. centered noise with unknown invertible
covariance matrix $\Gamma _{0}$.
\end{itemize}
Our goal is to estimate the true parameter by minimizing an appropriate
cost function. This model is called a regression model and a popular choice for the associated cost function is the mean square error : 
\[
\frac{1}{n}\sum _{t=1}^{n}\left\Vert Y_{t}-F_{W}\left(Z_{t}\right)\right\Vert ^{2}
\]
 where $\left\Vert .\right\Vert $ denotes the Euclidean norm on $\mathbb{R}^{d}$.
Although this function is widely used, it is easy to show that we get then a suboptimal estimator.
An other solution is to use an approximation of the covariance error matrix to compute generalized least square estimator :
\[
\frac{1}{n}\sum _{t=1}^{n}\left(Y_{t}-F_{W}\left(Z_{t}\right)\right)^{T}\Gamma ^{-1}\left(Y_{t}-F_{W}\left(Z_{t}\right)\right),
\]
 where  $T$ denotes the transposition of the matrix. Here we assume that $\Gamma $ is a good approximation of the true covariance matrix of the noise $\Gamma_{0}$.
However it takes time to compute a good approximation of matrix $\Gamma_{0}$ and it leads asymptotically to the cost function proposed in this article (see for example Rynkiewicz \cite{Rynkiewicz4}) : 
\begin{equation}
U_{n}\left(W\right):=\log \det \left(\frac{1}{n}\sum _{t=1}^{n}(Y_{t}-F_{W}(Z_{t}))(Y_{t}-F_{W}(Z_{t}))^{T}\right)
\end{equation}

This paper is devoted to the theoretical study of $U_{n}\left(W\right)$. We assume that the true architecture of the MLP is known so that the Hessian matrix computed in the sequel verifies the assumption to be definite positive  (see  Fukumizu \cite{Fukumizu}).

In this framework, we study the asymptotic behavior $\hat W_{n}:=\arg \min U_{n}\left(W\right)$,
the weights minimizing the cost function $U_{n}\left(W\right)$.
We show that under simple assumptions this estimator is asymptotically optimal in the sense that it has the same asymptotic behavior than the generalized least square estimator using  the true covariance matrix of the noise.

Numerical procedures to compute this estimator and examples of it behavior can be found in Rynkiewicz \cite{Rynkiewicz4}.
\section{The first and second derivatives of $W\longmapsto U_{n}\left(W\right)$}
First, we introduce a notation :  if $F_W(X)$ is a $d$-dimensional parametric function depending of a parameter $W$, let us write $\frac{\partial F_W(X)}{\partial W_k}$ (resp. $\frac{\partial^2 F_W(X)}{\partial W_k\partial W_l}$) for the $d$-dimensional vector of partial derivative (resp. second order partial derivatives) of each component of $F_W(X)$. 
\subsection{First derivatives}
Now, if $\Gamma_n(W)$ is a matrix depending of the parameter vector $W$, we get  From Magnus and Neudecker \cite{Magnus} \[
\frac{\partial }{\partial W_{k}}\ln \det \left(\Gamma _{n}(W)\right)=tr\left(\Gamma_{n}^{-1}(W) \frac{\partial }{\partial W_{k}}\Gamma _{n}(W)\right)\]
here 
\[
\Gamma_n(W)=\frac{1}{n}\sum _{t=1}^{n}(y_{t}-F_{W}(z_{t}))(y_{t}-F_{W}(z_{t}))^{T}
\] 
note that these matrix $\Gamma_{n}(W)$ and it inverse  are symmetric. 
Now, if we note 
\[
A_n(W_k)=\frac{1}{n}\sum _{t=1}^{n}\left(-\frac{\partial F_{W}(z_{t})}{\partial W_{k}}(y_t-F_{W}(z_{t}))^T\right)
\]
using the fact 
\[
tr\left(\Gamma_{n}^{-1}(W)A_n(W_k)\right)=tr\left(A_n^T(W_k)\Gamma_{n}^{-1}(W)\right)=tr\left(\Gamma_{n}^{-1}(W)A_n^T(W_k)\right)
\]

we get 
\begin{equation}\label{first_deriv}
\frac{\partial }{\partial W_{k}}\ln \det \left(\Gamma _{n}(W)\right)=2tr\left(\Gamma_{n}^{-1}(W)A_n(W_k)\right)
\end{equation}
\subsection{Second derivatives}
We write now 
\[
B_n(W_k,W_l):=\frac{1}{n}\sum _{t=1}^{n}\left( \frac{\partial F_{W}(z_{t})}{\partial W_{k}}\frac{\partial F_{W}(z_{t})}{\partial W_{l}}^T\right)
\]
and
\[
C_n(W_k,W_l):=\frac{1}{n}\sum _{t=1}^{n}\left( -(y_t-F_{W}(z_{t})) \frac{\partial^2 F_{W}(z_{t})}{\partial W_{k}\partial W_{l}}^T\right)
\]
We get

\[
\begin{array}{l}
\frac{\partial^2 U_n(W)}{\partial W_k\partial W_l}=\frac{\partial }{\partial W_{l}}2tr\left(\Gamma_{n}^{-1}(W)A_n(W_k)\right)=\\
2tr\left(\frac{\partial\Gamma_{n}^{-1}(W)}{\partial W_{l}}A(W_k)\right)+2tr\left(\Gamma_{n}^{-1}(W)B_n(W_k,W_l)\right)+2tr\left(\Gamma_{n}(W)^{-1}C_n(W_k,W_l)\right)\\
\end{array}
\]
Now, Magnus and Neudecker \cite{Magnus} give an analytic form of the derivative of an inverse matrix, so we get
\[
\begin{array}{l}
\frac{\partial^2 U_n(W)}{\partial W_k\partial W_l}=2tr\left(\Gamma_{n}^{-1}(W)\left(A_n(W_k)+A_n^T(W_k)\right)\Gamma_{n}^{-1}(W)A_n(W_k)\right)+\\
2tr\left(\Gamma_{n}^{-1}(W)B_n(W_k,W_l)\right)+2tr\left(\Gamma_{n}^{-1}(W)C_n(W_k,W_l)\right)
\end{array}
\]
so
\begin{equation}\label{second_deriv}
\begin{array}{l}
\frac{\partial^2 U_n(W)}{\partial W_k\partial W_l}=4tr\left(\Gamma_{n}^{-1}(W)A_n(W_k)\Gamma_{n}^{-1}(W)A_n(W_k)\right)\\
+2tr\left(\Gamma_{n}^{-1}(W)B_n(W_k,W_l)\right)+2tr\left(\Gamma_{n}^{-1}(W)C_n(W_k,W_l)\right)
\end{array}
\end{equation}
\section{Asymptotic properties of $\hat W_n$}
First, following the same lines that Yao \cite{Yao}, it is easy to show that, if the noise of the model has a moment of order at least 2, the estimator is strongly consistent (i.e. $\hat W_n \stackrel{a.s.}{\rightarrow} W^0$). 

Moreover, for a MLP function, there exists a constant $C$ such that we have the following inequalities :
\[
\begin{array}{l}
\Vert \frac{\partial F_W(Z)}{\partial W_k}\Vert\leq C(1+\Vert Z\Vert)\\
\Vert \frac{\partial^2 F_W(Z)}{\partial W_k\partial W_l}\Vert\leq C(1+\Vert Z\Vert^2)\\
\Vert \frac{\partial^2 F_W(Z)}{\partial W_k\partial W_l}-\frac{\partial^2 F_W^0(Z)}{\partial W_k\partial W_l}\Vert\leq C\Vert W-W^0 \Vert(1+\Vert Z\Vert^3)
\end{array}
\]
So, if $Z$ has a moment of order at least 3 (see the justification in Yao \cite{Yao}), we get the following lemma :

\begin{Lem}

Let $\Delta U_n(W^0)$ be the gradient vector of $U_n(W)$ at $W^0$,  $\Delta U(W^0)$ be the gradient vector of $U(W):=\log\det(Y-F_W(Z))$ at $W^0$ and $HU_n(W^0)$ be the Hessian matrix of  $U_n(W)$ at $W^0$. 

We define finally
\[
B(W_k,W_l):=\frac{\partial F_{W}(Z)}{\partial W_{k}}\frac{\partial F_{W}(Z)}{\partial W_{l}}^T
\]
and 
\[
A(W_k)=\left(-\frac{\partial F_{W}(Z)}{\partial W_{k}}(Y-F_{W}(Z))^T\right)
\]
We get then 
\begin{enumerate}
\item $HU_n(W^0)\stackrel{a.s.}{\rightarrow}2I_0$
\item $\sqrt{n}\Delta U_n(W^0)\stackrel{Law}{\rightarrow}{\cal N}(0,4I_0)$
\end{enumerate}
where, the component $(k,l)$ of the matrix $I_0$ is :
\[
tr\left(\Gamma^{-1}_0E\left(B(W^0_k,W^0_l)\right) \right)
\]
\end{Lem}
\paragraph{proof}
To prove the lemma, we remark first that the component $(k,l)$ of the matrix $4I_0$ is :
\[
E\left(\frac{\partial U(W^0)}{\partial W_k}\frac{\partial U(W^0)}{\partial W^0_l}\right)=E\left(2tr\left(\Gamma^{-1}_0A^T(W^0_k)\right)\times2tr\left(\Gamma^{-1}_0A(W^0_l)\right)\right)
\]
and, since the trace of the product is invariant by circular permutation, 
\[
\begin{array}{l}
E\left(\frac{\partial U(W^0)}{\partial W_k}\frac{\partial U(W^0)}{\partial W^0_l}\right)=\\
4E\left( -\frac{\partial F_{W^0}(Z)^T}{\partial W_k}\Gamma^{-1}_0(Y-F_{W^0}(Z))(Y-F_{W^0}(Z))^T\Gamma^{-1}_0\left(-\frac{\partial F_{W^0}(Z))}{\partial W_l}\right)\right)\\
=4E\left(\frac{\partial F_{W^0}(Z)^T}{\partial W_k}\Gamma^{-1}_0\frac{\partial F_{W^0}(Z)}{\partial W_l}\right)\\
=4tr\left(\Gamma^{-1}_0E\left(\frac{\partial F_{W^0}(Z)}{\partial W_k}\frac{\partial F_{W^0}(Z)^T}{\partial W_l}\right) \right)\\
=4tr\left(\Gamma^{-1}_0E\left(B(W^0_k,W^0_l)\right)\right) 

\end{array}
\]
Now, for the component $(k,l)$ of the expectation of the Hessian matrix, we remark that 
\[
\lim_{n\rightarrow \infty}tr\left(\Gamma_{n}^{-1}(W^0)A_n(W^0_k)\Gamma_{n}^{-1}(W^0)A_n(W^0_k)\right)=0
\]
and
\[
\lim_{n\rightarrow \infty}tr\Gamma_{n}^{-1}C_n(W^0_k,W^0_l)=0
\]
so
\[
\begin{array}{l}
\lim_{n\rightarrow \infty}H_n(W^0)=\lim_{n\rightarrow \infty}4tr\left(\Gamma_{n}^{-1}(W^0)A_n(W^0_k)\Gamma_{n}^{-1}(W^0)A_n(W^0_k)\right)+\\
2tr\Gamma_{n}^{-1}(W^0)B_n(W^0_k,W^0_l)+2tr\Gamma_{n}^{-1}C_n(W^0_k,W^0_l)=\\
=2tr\left(\Gamma^{-1}_0E\left(B(W^0_k,W^0_l)\right)\right)\\
\blacksquare
\end{array}
\]

From a classical argument of local asymptotic normality (see for example Yao \cite{Yao}), we deduce then the following property for the estimator $\hat W_n$ :  

\begin{Prop}
Let $W_n^*$ the estimator of the generalized least square :
\[
W_n^*:=\arg\min \frac{1}{n}\sum _{t=1}^{n}\left(Y_{t}-F_{W}\left(Z_{t}\right)\right)^{T}\Gamma_0^{-1}\left(Y_{t}-F_{W}\left(Z_{t}\right)\right)
\]
then we have
\[
\lim_{n\rightarrow\infty}\sqrt{n}(W_n^*-W^0)=\lim_{n\rightarrow\infty}\sqrt{n}(\hat W_n-W^0)={\cal N}(0,I_0^{-1})
\]

\end{Prop}

We remark that $\hat W_n$ has the same asymptotic behavior than the estimator  generalized
least square estimator with the true covariance matrix  $\Gamma^{-1}_{0}$ which is asymptotically optimal (see for example Ljung \cite{Ljung}),  so the proposed estimator is asymptotically optimal too.

\section{Conclusion}
In the linear multidimensional regression model the optimal estimator has an analytic solution (see Magnus and Neudecker \cite{Magnus}), so it doesn't make sense to consider minimization of a cost function. However, for the non-linear multidimensional regression model the ordinary least square estimator is sub-optimal if the covariance matrix of the noise is not the
identity matrix. We can overcome this difficulty by using the cost function $U_n(W)$. The numerical computation and the empirical properties of these estimator have been studied in a previous article (see rynkiewicz \cite{Rynkiewicz4}). In this paper, we have given a proof of  the optimality of the estimator associated with  $U_n(W)$. This is then a good choice for the estimation of multidimensional non-linear regression model with multilayer perceptron.

\begin{footnotesize}


\end{footnotesize}


\end{document}